\documentclass[10pt]{article}

\usepackage[T1]{fontenc}
\usepackage[latin1]{inputenc}
\usepackage[english]{babel}
\usepackage{graphicx}
\usepackage{amsmath,amscd}
\usepackage{amssymb}
\usepackage{amsfonts}
\usepackage{amsthm}
\usepackage{bbm}
\usepackage{geometry}
\usepackage{enumerate}
\usepackage{theoremref}
\usepackage{tikz-cd}  
\usepackage{ytableau}

\theoremstyle{plain}
\newtheorem{theo}{Theorem} \numberwithin{theo}{section} 
\newtheorem{lem}[theo]{Lemma}
\newtheorem{cor}[theo]{Corollary}
\newtheorem{propo}[theo]{Proposition}
\newtheorem{problem}[theo]{Problem}
\newtheorem{assumption}[theo]{Assumption}
\newtheorem{conjecture}[theo]{Conjecture}

\theoremstyle{definition}
\newtheorem{define}[theo]{Definition}

\theoremstyle{remark}
\newtheorem{bem}[theo]{Remark}
\newtheorem{bsp}[theo]{Example}

\newcommand{\definition}[1]{\begin{define}#1\end{define}}
\newcommand{\remark}[1]{\begin{bem}#1\end{bem}}
\newcommand{\example}[1]{\begin{bsp}#1\end{bsp}}
\newcommand{\pf}[1]{\begin{proof}#1\end{proof}}
\newcommand{\theorem}[1]{\begin{theo}#1\end{theo}}
\newcommand{\proposition}[1]{\begin{propo}#1\end{propo}}
\newcommand{\lemma}[1]{\begin{lem}#1\end{lem}}
\newcommand{\corollary}[1]{\begin{cor}#1\end{cor}}

\newcommand{\aufzaehlung}[1]{
\begin{description}
#1
\end{description}
}

\newcommand{\enumeration}[1]{
\begin{enumerate}
	#1
\end{enumerate}
}

\newcommand{\case}[1]{
\begin{cases}
#1
\end{cases}
}

\newcommand{\eqn}[1]{\begin{equation}#1\end{equation}}

\newcommand{\eqnlines}[2]{\begin{equation*}
\begin{array}{#1}#2
\end{array}
\end{equation*}
}

\newcommand{\C}{\mathbb{C}}   
\newcommand{\R}{\mathbb{R}}   
\newcommand{\Z}{\mathbb{Z}}   
\newcommand{\N}{\mathbb{N}}   
 
\newcommand{\SH}{\mathcal{H}} 
\newcommand{\SK}{\mathcal{K}} 
\newcommand{\SB}{\mathcal{B}} 
\newcommand{\SW}{\mathcal{W}} 

\newcommand{\LG}{\mathfrak{g}}   
\newcommand{\LH}{\mathfrak{h}}   
\newcommand{\LB}{\mathfrak{b}}   
\newcommand{\LU}{\mathfrak{u}}   
\newcommand{\LSO}{\mathfrak{so}} 
\newcommand{\LSU}{\mathfrak{su}} 
\newcommand{\LGL}{\mathfrak{gl}} 

\newcommand{\fL}{\textbf{L}} 

\DeclareMathOperator*{\I}{i}     

\newcommand{\eps}{\varepsilon}

\DeclareMathOperator*{\Tr}{Tr}     
\DeclareMathOperator*{\Hom}{Hom}   
\DeclareMathOperator*{\End}{End}   
\DeclareMathOperator*{\GL}{GL}     
\DeclareMathOperator*{\diag}{diag} 
\DeclareMathOperator*{\id}{id}     
\DeclareMathOperator*{\sgn}{sgn}   

\newcommand{\mynorm}[2]{\left\|#1\right\|_{#2}}
\newcommand{\opnorm}[1]{\mynorm{#1}{\text{op}}}
\newcommand{\skp}[3]{\langle #1,#2 \rangle_{#3}}
\newcommand{\gen}[2]{\left\langle #1\right\rangle_{#2}}

\newcommand{\indlim}{\varinjlim}

\newcommand{\hotimes}{\widehat{\otimes}}
\newcommand{\hoplus}{\widehat{\oplus}}
\newcommand{\oline}[1]{\overline{#1}}
\newcommand{\iof}{if and only if }                           
\newcommand{\1}{\mathbbm{1}}                                 

\newcommand{\cH}{\SH}
\newcommand{\urep}[1]{(\pi_{#1},\SH_{#1})}
\newcommand{\starurep}[1]{(\pi^*_{#1},\SH^*_{#1})}
\newcommand{\setlambda}{(\N_0)_\downarrow^{(\N)}}
\newcommand{\ZC}{Z^1_{\text{cond}}(G,\pi,\SH)}
\newcommand{\HC}{H^1_{\text{cond}}(G,\pi,\SH)}

\begin{document}

\title{On the First Order Cohomology of \\Infinite--Dimensional Unitary Groups} 

\author{Manuel Herbst and Karl-Hermann Neeb}

\maketitle

\begin{abstract}
The irreducible unitary highest weight representations $\urep{\lambda}$ of the group $U(\infty)$, which is the countable direct limit of the compact unitary groups $U(n)$, are classified by the orbits of the weights $\lambda \in \Z^{\N}$ under the Weyl group $S_{(\N)}$ of finite permutations. Here, we determine those weights $\lambda$ for which the first cohomology space $H^1(U(\infty),\pi_\lambda,\SH_\lambda)$ vanishes. For finitely supported $\lambda \neq 0$, we find that the first cohomology space $H^1(U(\infty),\pi_\lambda,\SH_\lambda)$ never vanishes. For these $\lambda$, the highest weight representations extend to norm-continuous irreducible representations of the full unitary group $U(\SH)$ (for $\SH:= \ell^2(\N,\C)$) endowed with the strong operator topology and to norm-continuous representations of the unitary groups $U_p(\SH)$ ($p\in [1,\infty]$) consisting of those unitary operators $g\in U(\SH)$ for which $g-\1$ is of $p$th Schatten class. However, not every 1-cocycle on $U(\infty)$ automatically extends to one on these unitary groups, so we may not conclude that the first cohomology spaces of the extended representations are non-vanishing. On the contrary, for the groups $U(\SH)$ and $U_\infty(\SH)$, all first cohomology spaces vanish. This is different for the groups $U_p(\SH)$ with $1\leq p <\infty$, where only the identical representation on $\SH$ and its dual representation have vanishing first cohomology spaces.

\end{abstract}

{\bf Key words:}
First order group cohomology, unitary representation, (Banach-) Lie group, Lie algebra, direct limit group 

\section{Introduction}
\label{sec: Intro}

For a continuous unitary or orthogonal representation $(\pi,\SH)$ of the topological group $G$ we call a map $\beta:G\to\SH$ a \textit{1-cocycle} if it satisfies the 1-cocycle relation \eqn{\label{1-cocycle} \beta(gh) = \beta(g) +\pi(g)\beta(h),\text{ for all } g,h\in G.} A 1-cocycle of the form $\beta(g)= \partial_v(g):= \pi(g)v-v$ for some $v\in \SH$ is called a \textit{1-coboundary} or a \textit{trivial} 1-cocycle. The vector space of all continuous 1-cocycles associated to $(\pi,\SH)$ is denoted by $Z^1(G,\pi,\SH)$, its subspace of 1-coboundaries by $B^1(G,\pi,\SH)$. The quotient space \[H^1(G,\pi,\SH):= Z^1(G,\pi,\SH)/B^1(G,\pi,\SH)\] is called the \textit{1-cohomology space}. \\

Interest in the first order cohomology of Lie groups has grown in the last decades since it is related to mathematical problems occurring in a wide range of mathematical disciplines such as geometric group theory, unitary representations, ergodic theory, stochastic processes and theoretical physics. Its motivation comes from the construction of irreducible unitary representations in Fock spaces (cf. \cite{Is96}, \cite{PS72}), the study of unitary representations of mapping groups/current groups (cf. \cite{Is96}, \cite{PS72}, \cite{Ar69}, \cite{A-T93}), the classification theory of negative definite functions (cf. \cite{FH74}), L\'evy processes and infinitely divisible probability distributions (cf. \cite{PS72}) and continuous tensor products of unitary representations (cf. \cite{PS72b}, \cite{St69}, \cite{Gui72}). The study of group cohomology spaces was notably propagated in the 1970s due to influential papers by J.P. Serre. One of the key results developed during that decade was the Delorme--Guichardet Theorem which states that, for a $\sigma$-compact locally compact group $G$, the following are equivalent (cf. \cite[Thm. 2.12.4, Prop. 2.2.10]{BHV08}):
\aufzaehlung{
\item[Property (T):] There exists a compact subset $K\subseteq G$ and some positive constant $\eps >0$ such that every continuous unitary representation $(\pi,\SH)$, for which there exists a unit vector $v\in\SH$ such that $\sup_{g\in K}\mynorm{\pi(g)v-v}{} < \eps$, has a non-trivial $G$-fixed vector, i.e. $\SH^G \neq \{0\}$.
\item[Property (FH):] $H^1(G,\pi,\SH) ={0}$ for every continuous orthogonal representation $(\pi,\SH)$ of $G$. 
}
More precisely, the implication $(T)\Longrightarrow (FH)$ holds for arbitrary topological groups whereas the $\sigma$-compactness is necessary for the converse implication (cf. \cite[Remark 2.12.5]{BHV08}).\\ 

It is a well-known fact that a 1-cocycle is trivial \iof it is bounded (cf. \cite[Prop. 2.2.9]{BHV08}). This implies that compact groups have property (FH) (see also \cite[Thm. 15.1]{PS72}). In the realm of irreducible unitary representations of finite dimensional Lie groups, the vanishing of the first cohomology spaces occurs surprisingly often: In view of the Delorme--Guichardet Theorem, \cite[Theorem 3.5.4]{BHV08} states that a connected semi-simple Lie group $G$ has property (FH) \iof no simple factor of its Lie algebra is isomorpic to $\LSU(n,1)$ or $\LSO(n,1)$ (for $n\in\N$). Moreover, up to equivalence, there exists only a finite number of (topologically complete) irreducible continuous unitary representations $\urep{}$ of a connected semi-simple Lie group, for which \\$H^1(G,\pi,\SH) \neq \{0\}$ (cf. \cite[Prop. 11]{PS75}). For a connected nilpotent Lie group, one can show that $H^1(G,\pi,\SH) =\{0\}$ holds for every nontrivial, irreducible continuous unitary representation $\urep{}$ (e.g. \cite[Thm. 17.4]{PS72}).\\

In this article, we turn to infinite dimensional Lie groups: Let $\SH\cong \ell^2(\N,\C)$ be an infinite dimensional complex separable Hilbert space. We consider the unitary groups 
\eqn{\label{rigging} U(\infty) \subset U_1(\SH) \subset\ldots \subset U_p(\SH) \subset\ldots \subset U_\infty(\SH) \subset U(\SH).}
Here $U(\infty)$ is the increasing union of the compact Lie groups $U(n):= U(n,\C)$, endowed with the direct limit Lie group structure (cf. \cite[Thm. 4.3]{Gloe05}). Its Lie algebra 
$\LU(\infty) = \indlim \LU(n)$ is an
increasing union of finite dimensional compact Lie algebras (cf. \cite[Thm. 4.3]{Gloe05}).
The other groups $U_p(\cH):= U (\cH)\cap(\1 + \SB_p(\cH))$ 
(for $1\leq p\leq\infty$) are Banach--Lie groups whose Lie algebras 
$\LU_p(\SH) = \LU(\cH) \cap \SB_p(\cH)$ are determined by the $p$th Schatten ideals $\SB_p(\SH)$ (cf. \cite[Sections II.5, II.6]{dlH72}). Note that $\SB_\infty(\cH)$ is the space of compact operators on $\cH$. In \cite{Ne98}, unitary highest weight representations of the group $U(\infty)$ have been classified by the orbits of the weights $\lambda = (\lambda_n)_{n \in \N}$ in $\Z^{\N}$ under the Weyl group $\SW \cong S_{(\N)}$. This means that two highest weight representations $\urep{\lambda}$ and $\urep{\mu}$ are equivalent \iof the entries of the weights $\lambda$ and $\mu$ coincide up to a finite permutation. The unitary representation $\urep{\lambda}$ extends to the group $U_1(\SH)$ \iof the weight $\lambda$ is bounded (cf. Proposition III.7 and Theorem III.4 in \cite{Ne98}). Moreover, the proof of Proposition III.10 in \cite{Ne98} shows that $\urep{\lambda}$ extends to the group $U_p(\SH)$ for some $p>1$ \iof the weight $\lambda$ has finite support. The problem is to determine for which weights $\lambda$, the corresponding first cohomology spaces are trivial. The case $\lambda =0$ corresponds to the one-dimensional trivial representation and the fact that there is no nonzero continuous group homomorphism $G\to\R$ for each of the unitary groups $G$ from \eqref{rigging} implies that the corresponding 1-cohomology spaces vanish (cf. \thref{sec: direct limit of compact groups}). For $\lambda\neq 0$, we derive the following solutions:
\aufzaehlung{
\item[Theorem 4.10:] For the direct limit $U(\infty)$, the first cohomology space $H^1(U(\infty),\pi_\lambda,\SH_\lambda)$ vanishes if and only if either all but finitely many entries of $\lambda$ are positive integers or all but finitely many entries are negative integers. In particular, for finitely supported weights $\lambda$, the corresponding 1-cohomology spaces never vanish.
\item[Theorem 7.7:] For $1\leq p <\infty$ and for finitely supported nonzero weights $\lambda$, the first cohomology space $H^1(U_p(\SH),\pi_\lambda,\SH_\lambda)$ vanishes in exactly two cases, namely either if $U_p(\SH)$ acts on $\SH$ via the identical representation or if $U_p(\SH)$ acts on $\SH^\ast$ via the dual representation. This corresponds to the case $\lambda = (1,0,0,0,\ldots)$ (natural action) resp. $\lambda = (-1,0,0,0,\ldots)$ (conatural action).
}
In particular, \thref{H^1 for U_p} reveals that the group $U_p(\SH)$ (with $1\leq p <\infty$) neither has property (FH) nor has property (T). For $p=2$, the fact that $U_2(\SH)$ does not have property (T) has already been shown in \cite{Pe17}. The question whether $U_2(\SH)$ has property (FH) is however stated there as an open problem (see \cite[3.5]{Pe17}) which is answered by \thref{H^1 for U_p}. \\

\cite[Prop. 6.5, Remark 6.8(a)]{At91} shows that the full unitary group $U(\SH)$ and the group $U_\infty(\SH)$ belong to the class of \textit{bounded} groups. We call a topological group $G$ \textit{bounded} if, for every identity neighborhood $U\subseteq G$, there exists an integer $n\in\N$ and a finite subset $F\subseteq G$ such that $G\subseteq U^nF$. \footnote{Using the terminology of \cite{At91}, this can be expressed as $G$ being bounded in the right uniformity of the topological group $G$. If $G$ is connected, one can always use $F=\{e\}$ which leads some authors to adopt this as a requirement in their definition.} That $U(\SH)$ and $U_{\infty}(\SH)$ are bounded groups is also shown in \cite{Ne13} with a proof relying on the spectral calculus of operators on Hilbert spaces. It is immediate from the definition that bounded groups have property (FH) \footnote{Indeed, if $\beta:G\to\SH$ is a continuous 1-cocycle associated to a continuous unitary representation $(\pi,\SH)$ of a bounded group $G$, then ~$U:=\{g\in G :\ \ \mynorm{\beta(g)}{} <1\}$ is an open identity neighborhood and $G\subset U^nF$ for some $n\in\N$ and some finite subset $F$ implies ~$\sup_{g\in G}\mynorm{\beta(g)}{} \leq n + \max_{f\in F}\mynorm{\beta(f)}{}< \infty$.}. Hence, if ~$G\in \{U_\infty(\SH),U(\SH)\}$ and $\urep{}$ is a continuous unitary representation of ~$G$, then $H^1(G,\pi,\SH) = \{0\}$ (see also \cite{Ro09}). The unitary groups $U(\infty)$ and $U_p(\SH)$ for $1\leq p <\infty$ are not bounded because they have nontrivial 1-cohomology spaces. \\

The core ideas for the proof of \thref{condition for SH_lambda=0} are the following: Any 1-cocycle of one of the unitary groups above has the property that its restriction to one of the compact subgroups $U(n)$ yields a 1-coboundary. This leads us to the notion of a \textit{conditional} 1-cocycle in Section \ref{sec: direct limit}. If $(G_j)_{j\in J}$ is a family of subgroups of the group $G$, then we call a 1-cocycle \textit{conditional} if its restriction to every subgroup $G_j$ yields a 1-coboundary. A special case arises if $G=\indlim G_n$ is the direct limit of an increasing sequence of subgroups $G_n$. The statement that every conditional 1-cocycle (w.r.t. the family of the subgroups $G_n$) is a 1-coboundary is then equivalent to the statement that every $G_n$-fixed vector is already fixed by the whole group $G$ provided that $n$ is sufficiently large (cf. \thref{H^1=0}). In Section \ref{sec: U(infty)}, we introduce the highest weight representations $\urep{\lambda}$ of $U(\infty)$ as a direct limit of highest weight representations of the subgroups $U(n)$ and we directly conclude that the representation $\urep{\lambda}$ admits a $U(n)$-fixed vector \iof one of the highest weight subrepresentations of $U(k)$ (with $k\geq n$) admits a $U(n)$-fixed vector. This observation allows us to reformulate the problem as a matter of branching from $U(k)$ to $U(n)$. Thus, \thref{condition for SH_lambda=0} is obtained by applying the classical branching law as stated in \thref{branching law}.\\

Our strategy to prove \thref{H^1 for U_p} can be outlined as follows: For finitely supported $\lambda$, we realize the extended highest weight representations $\urep{\lambda}$ of the group $U_p(\SH)$ (for $1\leq p<\infty$) as subrepresentations of finite tensor product representations that are built from the identical and its dual action of $U_p(\SH)$ on $\SH$ respectively on the dual $\SH^*$ (cf. Section \ref{sec: U_p}). For these finite tensor product representations (with at least two factors), it is easy to construct nontrivial 1-cocycles as a countably infinite sum of 1-coboundaries (cf. Section \ref{sec: finite tensor product}). Projecting these 1-cocycles onto the corresponding highest weight submodules gives nontrivial 1-cocycles for the highest weight representations $\urep{\lambda}$ (cf. \thref{H^1 not 0 for U_p}). That the first cohomology spaces for the identical representations on $\SH$ and the dual representation on $\SH^*$ (cf. Section \ref{sec: natural representation}) vanish is shown in \thref{sec: Lie algebra} and \thref{identical representation}.\\

In a subsequent paper, we plan to treat the case where $\lambda$ is bounded with infinite support and the corresponding highest weight representation $\urep{\lambda}$ of $U(\infty)$ extends to the group $U_1(\SH)$. The problem here is that it is much harder to realize the representation $\urep{\lambda}$ concretely than it is in the finitely supported case. If the first cohomology space vanishes for the group $U(\infty)$, then it also does for the group $U_1(\SH)$ since $U(\infty)$ is dense in $U_1(\SH)$. But this is not the case for every bounded weight with infinite support and this situation remains, at the present stage, an open problem.

\section*{Acknowledgements}
We would like to thank the referee for pointing out the references \cite{At91} and \cite{Ro09}.

\section{Notation, terminology and preliminaries}
\label{sec: basics}

Throughout this article, we use the symbol $\SH$ for a complex Hilbert space. We follow the convention that the scalar product $\skp{\cdot}{\cdot}{\SH}$ of a complex Hilbert space $\SH$ is linear in the first argument and antilinear in the second. If $V$ is any complex vector space, we write $\oline{V}$ for the complex vector space with the same addition operation but with scalar multiplication $(\lambda, v)\mapsto \oline{\lambda} v$ for $v\in V$ and $\lambda \in\C$. We use the symbol $G$ for a (topological) group and write $e$ for its neutral element. If not otherwise stated, we always assume that any topological group is a Hausdorff space. A \textit{continuous unitary representation} $(\pi,\SH)$ of $G$ on a complex Hilbert space $\SH$ is a group homomorphism $\pi:G\to U(\SH)$ which is continuous with respect to the strong operator topology on the unitary group $U(\SH)$. The unitary representation $(\pi,\SH)$ is called \textit{norm-continuous}, if $\pi$ is continuous with respect to the norm-topology on $U(\SH)$. We stick to the convention of writing $\SH^G$ for the subspace of $G$-fixed vectors.

\subsection{On the Lie groups used in this article}

Throughout this article, a Lie group $G$ is either a real Banach--Lie group or a direct limit of real finite-dimensional Lie groups. We write $\LG:= \fL(G) \cong T_e(G) $ for its (real) Lie algebra. The corresponding Lie bracket is denoted by $[\cdot,\cdot]$. Elements $z\in\LG_\C$ of the complexified Lie algebra split into $z= x+\I y$, where $x,y\in\LG$ are uniquely determined and are therefore called the \textit{real} resp. \textit{imaginary part} of $z$. We obtain an involution on $\LG_\C$ via $z^*:= -x+\I y$ for $z= x+\I y$.\\

Let $G$ be a Banach-Lie group and $(\pi,\SH)$ be a norm-continuous unitary representation of $G$. Every continuous 1-cocycle $\beta:G\to\SH$ defines a continuous group homomorphism of Banach-Lie groups 
\[\alpha: G\to \SH\rtimes U(\SH),\quad g\mapsto (\beta(g),\pi(g)).\]
By the Automatic Smoothness Theorem of Banach--Lie groups \footnote{For more details see \cite[Thm. IV.1.18]{Ne06} and the references given there.} we thus obtain a smooth map. In particular, both maps 
\[\pi:G\to U(\SH),\quad \beta:G\to\SH\] are smooth. Taking the derivative at $e$, we obtain continuous $\R$-linear maps
\[d\pi: \LG \to \LU(\SH),\quad d\beta:\LG \to \SH\]
satisfying, for arbitrary $x,y\in\LG$, the relations 
\[d\pi([x,y]) = d\pi(x)d\pi(y) - d\pi(y)d\pi(x),\quad d\beta([x,y]) = d\pi(x)d\beta(y)- d\pi(y)d\beta(x).\]
Complexification yields continuous $\C$-linear maps
\[d\pi_\C: \LG_\C \to \SB(\SH),\quad d\beta_\C:\LG_\C \to \SH.\]

\subsection{Unitary highest weight representations of the group $U(n)$}

Let $G$ be a compact Lie group. Standard references for the classification of irreducible unitary representations of $G$ in terms of highest weights are \cite{Ze73} (for the case $G=U(n)$) and \cite{GW98} (for the classical groups). For the infinite-dimensional unitary groups, we refer to the survey in \cite{Ne04b}, where the reader also finds an introduction to Banach--Lie groups and crucial results from \cite{Ne98}. Here, we briefly introduce the terms \textit{highest weight} and \textit{highest weight representation} for the group $U(n)$ because these concepts are mentioned several times throughout this paper.\\

Let $\urep{}$ be an irreducible, continuous unitary representation of $U(n)$. Then, $\SH$ is finite--dimensional and the derived Lie algebra representation $(d\pi_\C,\SH)$ of $\LGL(n,\C)$ is also irreducible. A \textit{weight} $\lambda$ of the representation $(d\pi_\C,\SH)$ is a linear functional $\lambda: \LH\to\C$ on the maximal abelian subalgebra $\LH$ of diagonal matrices of the Lie algebra $\LGL(n,\C)$. Let $\LB$ be a maximal solvable Lie subalgebra of $\LGL(n,\C)$ containing $\LH$. It is a consequence of Lie's Theorem on the finite dimensional representations of solvable Lie algebras that there exists a (unique) one-dimensional $\LB$-eigenspace in $\SH$. The corresponding $\LB$-eigenvalue is a Lie algebra homomorphism $\lambda:\LB\to\C$ and its restriction to $\LH$ is a weight which is called a \textit{highest weight (for $\LB$)} of the representation $(\pi,\SH)$. For instance, we may choose $\LB$ to be the upper triangular matrices in $\LGL(n,\C)$. Then, one can show that for $x= \diag(x_1,\ldots,x_n) \in \LH$, we have $\lambda(x) = \sum_{i=1}^n \lambda_i x_i$ with decreasingly ordered integer-valued coefficients $\lambda_i$. This means that any irreducible, continuous unitary representation of $U(n)$ defines a decreasingly ordered integer-valued $n$-tuple. Conversely, to every $\lambda\in\Z^n$ with $\lambda_1\geq \ldots\geq \lambda_n$ there exists an irreducible unitary representation $\urep{\lambda}$ of $U(n)$ with highest weight $\lambda$ (for the solvable subalgebra of upper triangular matrices). We call $\urep{\lambda}$ a \textit{unitary highest weight representation}. Two highest weight representations $\urep{\lambda}$, $\urep{\mu}$ are isomorphic \iof $\lambda = \mu$. Therefore, we obtain a one to one correspondence between the equivalence classes of irreducible, continuous unitary representations of $U(n)$ and the decreasingly ordered vectors in $\Z^n$. \\

For tuples $\lambda\in\Z^n$ which are not necessarily decreasingly ordered, we write $\urep{\lambda}$ for the unitary highest weight representation corresponding to the decreasingly ordered tuple $\lambda'$ which is obtained after a suitable permutation of the entries of $\lambda$. This has the following background: Every tuple $\lambda\in\Z^n$ occurs as a highest weight of some irreducible unitary representation of $U(n)$ for some solvable subalgebra $\LB$ containing $\LH$. It is a general fact in highest weight theory that two highest weight representations are isomorphic \iof the corresponding weights belong to the same orbit under the Weyl group action on the set of weights. The Weyl group action on $\Z^n$ can be realized by the natural permutation action of $S_n$ on $\Z^n$.\\
 
The polar map 
\[p: U(n)\times \LU(n) \to \GL(n,\C),\quad (g,x)\mapsto g\exp(\I x) \]
is a diffeomorphism and therefore, every continuous unitary representation $\urep{}$ of $U(n)$ extends to a holomorphic\footnote{ in the sense that the derivative is complex linear.} representation of $\GL(n,\C)$ via
$\pi_\C(g\exp(\I x)):= \pi(g)\exp{\I d\pi_\C(x)}$. This extension is unique (universal complexification) and, in view of Weyl's Unitarian Trick, there is a one-to-one correspondence between continuous, irreducible unitary representations of $U(n)$ and irreducible, holomorphic representations of $\GL(n,\C)$. In particular, the classification of the equivalence classes of holomorphic irreducible representations of $\GL(n,\C)$ in terms of highest weights is the same as for the group $U(n)$. In this paper, we will use occasionally some results from the highest weight theory of the group $\GL(n,\C)$ (especially from \cite{GW98}). The preceding remark shows that they may be applied equally well to the group $U(n)$.

\section{Conditional 1-cocycles}
\label{sec: direct limit}

In this section, we set the foundation for our analysis of the 1-cocycles of the unitary groups $U(\infty)$ and $U_p(\SH)$ (for $1\leq p <\infty$) in Sections \ref{sec: U(infty)} and \ref{sec: natural representation}. Assume that a topological group $G$ admits an increasing sequence of subgroups (such that the union is dense in $G$). If the restriction of a 1-cocycle on $G$ to each of these subgroups yields a 1-coboundary, then we call it a \textit{conditional} 1-cocycle. This concept is based on the observation that the restriction of any 1-cocycle on the unitary groups $U(\infty)$ and $U_p(\SH)$ (for $1\leq p <\infty$) to any of the compact groups $U(n)$ yields a trivial 1-cocycle. Any conditional 1-cocycle may be viewed as a limit of a sequence of 1-coboundaries and the question whether the coboundary property is prerserved under this limit can be translated into handy necessary and sufficient criteria (cf. \thref{coboundary criterion}). \thref{H^1=0} applies in particular to the group $U(\infty)$ and provides a solution to the problem to decide whether a 1-cohomology space of $U(\infty)$ is trivial.\\

Let $G$ be a topological group with an ascending sequence of topological subgroups $(G_n)_{n\in\N}$ such that $G_\infty := \bigcup_{n\in\N} G_n$ is a dense subgroup in $G$. 

\definition{
Let $\urep{}$ be a continuous unitary representation of $G$.
A 1-cocycle $\beta:G\to\SH$ is called \textit{conditional} if its restriction to every subgroup $G_n$ is a 1-coboundary. In particular, every 1-coboundary on $G$ is a conditional 1-cocycle. We write $\ZC$ for the vector space of conditional 1-cocycles and $\HC$ for the quotient $\ZC/B^1(G,\pi,\SH)$. 
}

We consider a conditional 1-cocycle $\beta \in \ZC$. Restricting $\beta$ to the subgroup $G_n$ yields a 1-coboundary by assumption, i.e. $\beta|_{G_n} = \partial_{v_n}$ for some $v_n \in \SH$. The vector $v_n$ is unique up to adding a $G_n$-fixed vector. We shall write $\SH_n := \SH^{G_n}$ for the subspace of $G_n$-fixed vectors. Any conditional 1-cocycle $\beta$ thus defines a unique sequence $(v_n \in \SH_n^\bot)_{n\in\N}$ with the following compatibility condition
\eqn{\label{compatibility} v_n-v_m \in \SH_m \text{ if } m\leq n.} This means that we can write $\beta(g) = \lim_{n\rightarrow\infty} \pi(g)v_n-v_n$ for any $g\in G_\infty$. We will reformulate the compatibility condition in two ways:\\

First, note that the spaces $\SH_n$ form a decreasing sequence of Hilbert subspaces and the intersection is given by $ \bigcap_{n\in\N}\SH_n = \SH^G$, the subspace of $G$-fixed vectors. For $n=0$ we put $G_0 := \{e\}$, $\SH_0:=\SH$ and $v_0 :=0$. In view of $\SH_n^\bot \cong \oplus_{k=1}^n\left(\SH_k^\bot\cap\SH_{k-1}\right)$, we can write $v_n = \sum_{k=1}^n w_k$ for uniquely determined $w_k\in W_k:=\SH_k^\bot\cap\SH_{k-1}$. The compatibility condition \eqref{compatibility} then reads $w_k = v_k-v_{k-1}$, hence every conditional 1-cocycle $\beta$ defines a unique sequence $(w_n\in W_n)_{n\in\N}$. For $g\in G_\infty$, we then have \eqn{\label{w_n and beta} \beta(g) = \sum_{n\in\N}\pi(g)w_n-w_n.} 
A second way to interpret \eqref{compatibility} is as follows: For any $v\in \bigcup_{n\in\N} \SH_n^\bot$, the limit $\lim_{n\rightarrow \infty}\skp{v_n}{v}{}$ exists and we obtain an antilinear functional $a:\bigcup_{n\in\N} \SH_n^\bot \to \C,\quad v\mapsto \lim_{n\rightarrow\infty} \skp{v_n}{v}{}$ which has the property that its restriction to every subspace $\SH_n^\bot$ is continuous. Moreover, if $g\in G_\infty$ and $v\in\SH$, then $\pi(g)v-v\in \bigcup_{n\in\N}\SH_n^\bot$ and we obtain the relation
\eqn{\label{a and beta} a(\pi(g)v-v) = \lim_{n\rightarrow\infty} \skp{v_n}{\pi(g)v-v}{} = \lim_{n\rightarrow\infty} \skp{\pi(g^{-1})v_n-v_n}{v}{} = \skp{\beta(g^{-1})}{v}{}. }

\remark{\thlabel{sec: Lie algebra}
In a Lie-group context, we assume that $G$ is a Banach--Lie group with an increasing sequence $(G_n)_{n\in\N}$ of Banach-Lie subgroups such that $G_\infty = \bigcup_{n\in\N} G_n$ is a dense subgroup. The corresponding Banach-Lie algebras $\LG_n$ form an increasing sequence of Lie subalgebras of $\LG$. We write $\LG_\infty:= \bigcup_{n\in\N}\LG_n \subseteq \LG$ for the union which is a Lie subalgebra of $\LG$. 
Now, let $(\pi,\SH)$ be a norm-continuous unitary representation of $G$ and let $\beta \in \ZC$. Recall that we obtain continuous linear maps $d\pi_\C: \LG_\C \to \SB(\SH)$ and $d\beta_\C:\LG_\C\to\SH$ via complex linear extension of the derivatives. If $(v_n)_{n\in\N}$ is the sequence defined by $\beta$, then the compatibility condition \eqref{compatibility} implies $d\beta_\C(z) = \lim_{n\rightarrow\infty} d\pi_\C(z)v_n$. For any $z\in \LG_{\infty,\C}\subseteq \LG_\C$ and $v\in\SH$, we have $d\pi_\C(z)v\in \bigcup_{n\in\N}\SH_n^\bot$ and obtain
\eqn{\label{a and dbeta_C} a(d\pi_\C(z)v) = \lim_{n\rightarrow\infty} \skp{v_n}{d\pi_\C(z)v}{} = \lim_{n\rightarrow\infty} \skp{d\pi_\C(z^*)v_n}{v}{} = \skp{d\beta_\C(z^*)}{v}{}. }
This fact will be used later in Section \ref{sec: natural representation}.
}

\lemma{\thlabel{coboundary criterion}
Let $\beta$, $v_n$ and $w_n$ be as above.
For a vector $v\in \left(\SH^G\right)^\bot$, the following are equivalent:
\aufzaehlung{
\item[$i)$] The conditional 1-cocycle $\beta$ is a 1-coboundary with $\beta = \partial_v$.
\item[$ii)$] The sequence $(w_n)_{n\in\N}$ is square-summable and $v = \sum_{n\in\N}w_n$.
\item[$iii)$] The sequence $(v_n)_{n\in\N}$ converges in $\SH$ with limit $v$.
\item[$iv)$] The antilinear functional $a$ extends to a continuous antilinear functional on $\SH$ such that \\$a(v') = \skp{v}{v'}{}$ for all $v'\in \SH$.
}}

\pf{
\aufzaehlung{
\item[$i)\Longrightarrow ii)$]
For each vector $v_n = \sum_{k=1}^n w_k$, we have $v-v_n \in \SH_n = \SH^{G_n}$, hence $(v-v_n) \bot v_n$. Since the vectors $w_k$ are mutually orthogonal, we find $\mynorm{v}{}^2 \geq \sum_{k=1}^n \mynorm{w_k}{}^2$ for all $n\in\N$. This shows that the sequence $(w_k)_{k\in\N}$ is square integrable and the series $\sum_{k=1}^\infty w_k$ converges in $\SH$. By virtue of \eqref{w_n and beta}, we have $[\pi(g)-\1]v = [\pi(g)-\1]\sum_{k=1}^\infty w_k$ for each $g\in G$ and, since both $v$ and $\sum_{k=1}^\infty w_k$ belong to $\left(\SH^G\right)^\bot$, we conclude $v=\sum_{k=1}^\infty w_k$.
\item[$ii)\Longrightarrow iii)$] This is clear, since the convergence of the series implies that the sequence $(v_n)_{n\in\N}$ is a Cauchy-sequence in $\SH$ whose limit coincides with the limit of the series.
\item[$iii)\Longrightarrow iv)$] This follows directly from the definition of $a$.
\item[$iv)\Longrightarrow i)$] For any $v'\in\SH$ and $g\in G$, we have $\skp{\pi(g^{-1})v-v}{v'}{} = a(\pi(g)v'-v') = \skp{\beta(g^{-1})}{v'}{}$ by equation \eqref{a and beta}. This shows that $\beta = \partial_v$.
}}

\thref{coboundary criterion} shows that if all but finitely many of the spaces $W_n$ are trivial, then every conditional 1-cocycle is trivial. This occurs in particular, if the Hilbert space $\SH$ is finite dimensional. 
\remark{\thlabel{W_n and H^G_n}
We have $\left(\SH^{G_n}\right)^\bot = \oplus_{k=1}^n W_k$ for any $n\in\N$ and $\left(\SH^G\right)^\bot = \hoplus_{k\in\N} W_k$. This shows that, for any $n\in\N$, \[(\forall k>n)\quad W_k =\{0\} \Longleftrightarrow \SH^{G_n} = \SH^G. \]
}
Thus, we have shown the following corollary.
\corollary{\thlabel{H^G=H^G_n sufficient}
If $\SH^G = \SH^{G_n}$ for some $n\in\N$, then $\HC = \{0\}$.
}

We conclude our discussion on conditional 1-cocycles by specializing to the case where \\$G=\bigcup_{n\in\N}G_n$. We write $G=\indlim G_n$ whenever the group topology on $G$ coincides with the direct limit topology. This requirement ensures that every sequence $(w_n\in W_n)_{n\in\N}$ defines a conditional 1-cocycle via equation \eqref{w_n and beta}. Thus, the space $\ZC$ may be described in terms of the spaces $W_n$ and we conclude that the converse statement of \thref{H^G=H^G_n sufficient} holds. 

\proposition{\thlabel{H^1=0}
For $G= \indlim G_n$, we have the following statements
\enumeration{[i)]
\item $\ZC \cong \prod_{n\in\N}\left(\SH^{G_n}\right)^\bot\cap\SH^{G_{n-1}}$.
\item $\HC = \{0\} \Longleftrightarrow (\exists n\in\N)\quad \SH^{G_n} = \SH^G.$
}}

\pf{
\enumeration{[{a}d i)]
\item Every sequence $(w_n\in W_n)_{n\in\N}$ defines a map $\beta:G\to \SH$ via \\$\beta(g) := \sum_{n\in\N}\pi(g)w_n-w_n$ for $g\in G$. Then, the restriction to every subgroup $G_n$ yields a 1-coboundary $\beta|_{G_n} = \partial_{v_n}$ with $v_n = \sum_{k=1}^n w_k$. In particular, the restriction of $\beta$ to $G_n$ is continuous. Therefore, $\beta$ is continuous with respect to the direct limit topology on $G$. This shows that $\beta$ is a conditional 1-cocycle. The above mapping defines an isomorphism of vector spaces $\prod_{n\in\N}W_n \to \ZC$.
\item In view of \thref{H^G=H^G_n sufficient} and \thref{W_n and H^G_n}, it remains to show the implication \\$\HC = \{0\} \Longrightarrow (\exists n\in\N)(\forall k>n)\quad W_k = \{0\}$. The requirement $\ZC = B^1(G,\pi,\SH)$ amounts to saying that every sequence $(w_n)_{n\in\N}\in \prod_{n\in\N}W_n$ has to be square summable (\thref{coboundary criterion} and $i)$). This possible only if all but finitely many $W_n$ are trivial.\qedhere
}}



\remark{\thlabel{sec: direct limit of compact groups}
Assume that all subgroups $G_n$ are compact. Then, the restriction of every 1-cocycle on $G$ to one of the $G_n$ is a 1-coboundary, hence every 1-cocycle on $G$ is conditional. This means that $Z^1(G,\pi,\SH) = \ZC$ and $H^1(G,\pi,\SH) = \HC$. Note, that \thref{H^G=H^G_n sufficient} shows that for $\pi=\1$, one has $\{0\}= H^1(G,\1,\SH) \cong \Hom_{Grp}(G,\SH)$. This is a version of the fact that there is no nontrivial, continuous group homomorphism of $G$ into the group of additive real numbers. Moreover, the direct limit topology on $G_\infty = \bigcup_{n\in\N}G_n$, i.e. the finest group topology on $G_\infty$ for which all inclusions $G_n \hookrightarrow G_\infty$ are continuous, is a group topology (cf. \cite[Theorem 2]{Yam98}). Hence, \thref{H^1=0} applies in particular to the group $G_\infty$. The unitary group $U(\infty)$ is a direct limit of compact subgroups which are given by the unitary $n\times n$-matrices $U(n)$. Every unitary group $U_p(\SH)$ (for $\SH=\ell^2(\N,\C)$ and $1\leq p <\infty$) contains $U(\infty)$ as a dense subgroup. The one-dimensional trivial representation of the unitary groups $U(\infty)$ and $U_p(\SH)$ is the unitary highest weight representation that corresponds to the weight $\lambda =0$ and we have $H^1(G,\pi_0,\SH_0) = \{0\}$ for every unitary group $G$ occuring in \eqref{rigging}. 
}

\section{Unitary highest weight representations of $U(\infty)$}
\label{sec: U(infty)}

In this section, we define for every $\lambda\in\Z^\N$ an irreducible unitary representation $\urep{\lambda}$ of $U(\infty)$ that we call \textit{unitary highest weight representation (with highest weight $\lambda$)}. We determine for which $\lambda$ the spaces $H^1(U(\infty),\pi_\lambda,\SH_\lambda)$ are trivial (cf. \thref{condition for SH_lambda=0}). The main ingredients are the classical Branching Law for the highest weight representations of the unitary groups $U(n)$ (cf. \thref{branching law}) and \thref{H^1=0}.\\

Let $n\in\N$. Recall from Section \ref{sec: basics} that for any integer-valued tuple $\lambda\in\Z^n$, one can associate a unitary highest weight representation $\urep{\lambda}$ of $U(n)$. This representation is irreducible and every irreducible unitary representation of $U(n)$ is isomorphic to some $\urep{\lambda}$. Two highest weight representations $\urep{\lambda}$, $\urep{\mu}$ are isomorphic \iof the entries of the weights $\lambda$ and $\mu$ coincide up to permutation. For $n\geq 2$ the subgroup $U(n-1)$ decomposes a highest weight representation $\urep{\lambda}$ into a finite sum of highest weight representations $\urep{\eta}$ of $U(n-1)$ and it is natural to ask which weights $\eta\in \Z^{n-1}$ occur in the decomposition. The answer is a classical result in branching theory.

\theorem{\thlabel{branching law}
Let $\eta\in\Z^n$ and $\lambda\in\Z^{n+1}$ be two decreasingly ordered integer valued tuples. Then, the unitary highest weight representation $\urep{\eta}$ of $U(n)$ is a subrepresentation of $\urep{\lambda}$  \iof the tuples $\eta$ and $\lambda$ satisfy the interlacing condition
\eqn{\label{interlacing} \lambda_1\geq\eta_1\geq \lambda_2\geq\eta_2\geq\ldots\geq\lambda_n\geq \eta_n\geq\lambda_{n+1}.}
}
\pf{
A proof can be found e.g. in \cite[Thm 8.1.1]{GW98}.
}

\definition{
Let $\eta\in\Z^n$ and $\lambda\in\Z^{n+1}$ be two integer valued tuples. We say that $\eta$ \textit{interlaces} $\lambda$ and write $\eta\preccurlyeq \lambda$ if condition \eqref{interlacing} is satisfied after a suitable permutation of the entries of both tuples. 
}

Now, let $\lambda\in\Z^\N$ be an integer-valued sequence. For any $n\in\N$, we define the tuple $\lambda^{(n)}:= (\lambda_1,\lambda_2,\ldots,\lambda_n)\in\Z^n$ to be the $n$-tuple consisting of the first $n$ entries of $\lambda$. Then, $\lambda^{(n)}$ interlaces $\lambda^{(n+1)}$ and, according to \thref{branching law}, $\urep{\lambda^{(n)}}$ occurs as a sub-representation of $\urep{\lambda^{(n+1)}}$ (w.r.t. the group $U(n)$). Hence, each $\SH_{\lambda^{(n)}}$ is isometrically embedded into $\SH_{\lambda^{(n+1)}}$ and we write $V_\lambda:= \indlim \SH_{\lambda^{(n)}}$ for the direct limit in the category of pre-Hilbert spaces. The representations $\pi_{\lambda^{(n)}}$ canonically define a unitary representation of $U(\infty)= \bigcup_{n\in\N}U(n)$ on the pre-Hilbert space $V_\lambda$ which we denote by $\pi_\lambda$. It is clear that $\pi_\lambda$ extends to a unitary representation on the Hilbert space completion $\SH_\lambda:= \oline{V_\lambda}$ which is again denoted by $\pi_\lambda$.

\lemma{\thlabel{continuous, irreducible}
The unitary representation $\urep{\lambda}$ is continuous and irreducible. 
} 
\pf{
To see the continuity, it is enough to verify that, for any $v\in V_\lambda \subseteq \SH_\lambda$, the orbit map $g\mapsto \pi_\lambda(g)v$ is continuous. For fixed $v\in V_\lambda$ and sufficiently large $n$, we may assume that $v\in\SH_{\lambda^{(n)}}$. Therefore, the restriction of the orbit map to the subgroup $U(n)$ is continuous since the representation $\urep{\lambda^{(n)}}$ is continuous. This is true for all sufficiently large $n$ and the continuity thus follows from the fact that the group topology on $U(\infty)$ is given by the direct limit topology. That $\urep{\lambda}$ is irreducible follows from the fact that direct limits of irreducible representations are irreducible (cf. Proposition A.5 in \cite{BN12}).
}

\definition{\thlabel{HWR for U(infty)}
The unitary representation $\urep{\lambda}$ is called \textit{unitary highest weight representation} of $U(\infty)$ \textit{with highest weight $\lambda\in\Z^\N$}.
} 

\remark{
In \cite[Sections I,II]{Ne98}, the unitary highest weight representations of the direct limit Lie algebra $\LGL(\infty) = \bigcup_{n\in\N}\LGL(n)$ are classified in terms of real-valued sequences $\lambda\in \R^\N$. The corresponding highest weight module is denoted by $L(\lambda)$. In Section III of \cite{Ne98} it is shown that the underlying Lie algebra representation of $L(\lambda)$ integrates to a representation $\hat{\varrho}_\lambda: \GL(\infty)\to\End(L(\lambda))$ \iof $\lambda\in\Z^\N$. For any such integer-valued weight $\lambda$, the module $L(\lambda)$ may be identified with our pre-Hilbert space $V_\lambda$ and the restriction of $\hat{\varrho}_\lambda$ to the subgroup $U(\infty)$ is just our $\pi_\lambda$. Theorem I.20 of \cite{Ne98} states that two highest weight modules $L(\lambda)$ and $L(\mu)$ are equivalent \iof the weights belong to the same orbit under the Weyl group $\SW$. As remarked in Section~ II, the Weyl group for $\LGL(\infty)$ may be identified with the group $S_{(\N)}$ of finite permutations on the entries of the weights. Thus, we conclude that two unitary highest weight representations $\urep{\lambda}$ and $\urep{\mu}$ of $U(\infty)$ are equivalent \iof the entries of the weights $\lambda$ and $\mu$ coincide up to a finite(!) permutation of the entries.}

Using \thref{H^1=0}, we want to determine for which tuples $\lambda\in\Z^\N$, the 1-cohomology space for the corresponding highest weight representation $\urep{\lambda}$ vanishes. This is true for $\lambda=0$ (cf. \thref{sec: direct limit of compact groups}) and from now on, we assume that $\lambda\neq 0$. Note that then $\SH^{U(\infty)}_\lambda = \{0\}$ by the irreducibility of the highest weight representation.\\

We need the following two observations that follow directly from \thref{subrepresentation} below:
Let $\lambda\in\Z^\N$ and $\urep{\lambda}$ be the corresponding unitary highest weight representation. For $n\in\N$, we have 
\eqn{\label{From H_lambda to H_lambda_k}\SH^{U(n)}_\lambda = \{0\} \Longleftrightarrow (\forall k\geq n)\quad \SH^{U(n)}_{\lambda^{(k)}} = \{0\}.}
For $k>n$, we have 
\eqn{\label{From k>n to n+1>n} \SH^{U(n)}_{\lambda^{(k)}} \neq \{0\} \Longleftrightarrow \SH^{U(n)}_{\eta^{(k-1)}} \neq\{0\} \text{ for some }\eta^{(k-1)}\preccurlyeq \lambda^{(k)}.}
Both observations follow from the following lemma:
\lemma{\thlabel{subrepresentation}
Let $(\pi,\SH)$ be a continuous unitary representation of the topological group $G$. Assume that, for a subset $J\subseteq\N$, we have a family of $G$-invariant subspaces $(\SH_j)_{j\in J}$ such that one of the following conditions is satisfied:
\enumeration{[i)]
\item $\SH= \widehat{\oplus}_{j\in J}\SH_j$
\item $\SH = \oline{\bigcup_{j\in J}\SH_j}$
}
Let $(\rho,\SK)$ be an irreducible representation of $G$. Then, $(\rho,\SK)$ occurs as a subrepresentation in $(\pi,\SH)$ \iof it occurs as a subrepresentation in $(\pi,\SH_j)$ for some $j\in J$. 
}

\pf{
We denote by $P_j$ the orthogonal projection onto the subspace $\SH_j$. We may assume that $\SK\subset \SH$ is a Hilbert subspace of $\SH$. If $(\rho,\SK) \subset (\pi,\SH)$, then, in both cases, we find an index $j\in J$ for which $P_j(\SK)\neq \{0\}$. The projection operator $P_j|_{\SK}:\SK\to\SH_j$ intertwines $\rho$ and $\pi$. Therefore, $P_j(\SK) \subset \SH_j$ is a $G$-invariant subspace. By Schur's Lemma, we further conclude that $(P_j|_{\SK})^*P_j|_{\SK} = c\1|_{\SK}$. Since the operator $(P_j|_{\SK})^*P_j|_{\SK}$ is positive, the constant $c$ is real-valued and nonngeative. The case $c=0$ is excluded by $P_j(\SK)\neq \{0\}$. Therefore, $A:= \frac{1}{\sqrt{c}}P_j|_{\SK}:\SK\to P_j(\SK)$ is a linear isometry intertwining $\rho$ and $\pi$. This shows $(\rho,\SK)\subset (\pi,\SH_j)$. The converse statement is trival.
}

The equivalence \eqref{From k>n to n+1>n} motivates the following definition.
\definition{
Assume that a property (P) is defined for all tuples $\lambda \in \Z^k$ of length $k\geq n$. We say that (P) is \textit{interlacing-inheritable} if a tuple $\lambda$ has property (P) \iof there exists a tuple $\eta\preccurlyeq \lambda$ with property (P). 
}
The property that the highest weight representation $\urep{\lambda}$ of $U(k)$ corresponding to a tuple $\lambda\in\Z^k$ has a nontrivial $U(n)$-fixed vector (where $k\geq n$), is interlacing-inheritable. Any two interlacing-inheritable properties (P) and (P') are equivalent \iof they are equivalent for all tuples $\lambda\in \Z^n$ of length $n$.

\lemma{\thlabel{interlacing-inheritable}
Let $k\geq n$. For any $\lambda\in\Z^k$, the property 
\eqn{\label{property P} \#\{j: \lambda_j\geq 0\} \geq n \quad\text{ and }\quad \#\{j: \lambda_j\leq 0\} \geq n }
is interlacing-inheritable.
}

\pf{
Let $\lambda\in\Z^k$ for $k>n$. We may assume w.l.o.g. that $\lambda$ is decreasingly ordered. The interlacing condition $\eta\preccurlyeq\lambda$ then reads
$\lambda_1\geq\eta_1\geq\lambda_2\geq\eta_2\geq\ldots\geq\eta_{k-1}\geq\lambda_k$. If $\eta$ has the property \eqref{property P}, then the first (last) $n$ entries of $\eta$ are $\geq 0$ ($\leq 0$), hence so are the first (last) $n$ entries of $\lambda$. Conversely, assume that $\lambda$ has the required property. We construct $\eta\in\Z^{k-1}$ as follows: If $k>2n$, choose $\eta\preccurlyeq\lambda$ such that the first (last) $n$ entries of $\eta$ coincide with the first (last) $n$ entries of $\lambda$. If $k\leq 2n$, consider the set $J:= \{1,2,\ldots,n\}\cap\{k-1,k-2,\ldots,k-n\}$. Put $\eta_j:= \lambda_j$ whenever the index $j$ is smaller than the indices from $J$ and put $\eta_j:= \lambda_{j+1}$ whenever the index $j$ is greater than the indices from $J$.
}

\proposition{\thlabel{SH^G_n_lambda^k ungleich 0}
Let $k\geq n$ and $\lambda\in\Z^k$. Then the unitary highest weight representation $\urep{\lambda}$ (of $U(k)$) admits nonzero $U(n)$-fixed vectors \iof the following condition is satisfied:
\eqn{\label{weight property}\#\{j| \lambda_j \geq 0\} \geq n \quad\text{ and }\quad \#\{j| \lambda_j \leq 0\}\geq n. }
}
\pf{
In view of \thref{interlacing-inheritable}, both conditions are interlacing-inheritable and it just remains to check that they are equivalent for tuples $\lambda\in\Z^n$. Indeed, \eqref{weight property} is satisfied \iof $\lambda = 0$.
}

\theorem{\thlabel{condition for SH_lambda=0}
Let $\lambda\in\Z^{\N}\backslash\{0\}$ and $\urep{\lambda}$ be the corresponding unitary highest weight representation of $U(\infty)$. Then $H^1(U(\infty),\pi_\lambda,\SH_\lambda) = \{0\}$ \iof either \[ \#\{j| \lambda_j \leq 0\} < \infty \quad\text{ or }\quad \#\{j| \lambda_j \geq 0\} < \infty. \]
}

\pf{
According to \thref{H^1=0}, the first cohomology space vanishes \iof there exists some $n\in\N$ for which $\SH^{U(n)}_\lambda = \{0\}$. This is equivalent to $\SH^{U(n)}_{\lambda^{(k)}} = \{0\}$ for all $k\geq n$ (cf. \eqref{From H_lambda to H_lambda_k}), where $\lambda^{(k)}\in\Z^k$ consists of the first $k$ entries of $\lambda$. By virtue of \thref{SH^G_n_lambda^k ungleich 0} this can be rewritten as 
\eqnlines{lr}{
& (\forall k\geq n)\qquad \#\{j| \lambda^{(k)}_j \geq 0\} < n \quad\text{ or }\quad \#\{j| \lambda^{(k)}_j \leq 0\} < n. \\
\Longleftrightarrow &  \#\{j| \lambda_j \geq 0\} < n \quad\text{ or }\quad \#\{j| \lambda_j \leq 0\} < n. 
}}

\section{The identical representation of $U_p(\SH)$ on $\SH$}
\label{sec: natural representation}

Now, we prove that the natural actions of the group $U_p(\SH)$ on $\SH$ resp. on its topological dual $\SH^*$ have trivial 1-cohomology spaces. \\

Let $1\leq p<\infty$ and $\SH$ be a complex separable Hilbert space with ONB $(e_n)_{n\in\N}$. The unitary group $U_p(\SH) = U(\SH)\cap(\1 + \SB_p(\SH))$ is a Banach-Lie group with Banach-Lie algebra $\LU_p(\SH)=\LU(\SH)\cap \SB_p(\SH)$, i.e. the skew hermitian operators of $p$-th Schatten class (cf. \cite{dlH72}). Its complexification is given by the $p$-th Schatten operators $\LU_{p,\C}(\SH) = \SB_p(\SH)$. For any $z\in \LU_{p,\C}(\SH)$ the involution $z^*$ coincides with the usual Hilbert adjoint. The inductive limit $U(\infty)$ is a dense subgroup in $U_p(\SH)$ so that we can apply the results of Section \ref{sec: direct limit} and in particular \thref{sec: Lie algebra} to any norm-continuous unitary representation of $U_p(\SH)$. The complexification of the Lie algebra of $U(\infty)$ is given by $\LGL(\infty):= \indlim \LGL(n,\C)$. \\

The identical action of $U_p(\SH)$ on $\SH$ is given by the prescription \[U_p(\SH)\times \SH\to \SH,\quad (g,v)\mapsto gv.\] This action, restricted to $U(\infty)$, defines a unitary highest weight representation corresponding to the tuple $\lambda:= (1,0,0,\ldots)$ which we denote by $\urep{\lambda} = (\id,\SH)$. It is clear that the identical action defines a norm-continuous unitary representation of the group $U_p(\SH)$.\\

\proposition{\thlabel{identical representation}
We have $H^1(U_p(\SH),\id,\SH) = \{0\}$.
}

\pf{
For $u,w\in\SH$, we consider the rank-1-operator $w\otimes u^*:= \skp{\cdot}{u}{}w$ which is an element of $\SB_p(\SH)$ for each $p\in [1,\infty]$ since \[\mynorm{w\otimes u^*}{p} = \mynorm{w}{}\mynorm{u}{}.\]
For $v\in \SH_0 := \gen{e_n: n\in\N}{lin}$, we have $v\otimes e_1^* \in \LGL(\infty)$. Now, let $\beta \in Z^1(U_p(\SH),\id,\SH)$ and $a:= a_\beta$ be the corresponding antilinear functional from Section \ref{sec: direct limit}. Using equation \eqref{a and dbeta_C}, we find
\[|a(v)| = |a(v\otimes e_1^* (e_1))| = |\skp{d\beta_\C(e_1\otimes v^*)}{e_1}{}| \leq \mynorm{d\beta_\C}{}\mynorm{e_1\otimes v^*}{p} = \mynorm{d\beta_\C}{}\mynorm{v}{}\]
which shows that $a$ extends to a continuous linear functional on $\SH$. The assertion now follows from \thref{coboundary criterion}. 
}

The dual representation $(\id^*,\SH^*)$ is also a unitary highest weight representation which corresponds to the tuple $\lambda = (-1,0,0,\ldots)$.

\proposition{\thlabel{dual representation}
We have $H^1(U_p(\SH),\id^*,\SH^*) = \{0\}$.
}

\pf{
This immediately follows from \thref{identical representation} and the fact that the first order cohomology space for the dual representation vanishes \iof it  vanishes for the original representation.
}

\section{A simple way to construct unbounded 1-cocycles in finite tensor products}
\label{sec: finite tensor product}

This section should be viewed as a preparation for Section \ref{sec: U_p}. We consider a countably infinite sum of 1-coboundaries that converges pointwise and ask whether we thus obtain an unbounded 1-cocycle. To make life easier, we assume the underlying group $G$ to be completely metrizable since then, the pointwise converging sum is automatically continuous. This observation is based on a Baire category argument which is carried out in Appendix \ref{sec: infinite sums}. Note that the unitary groups $U_p(\SH)$ are completely metrizable for any $p\in [1,\infty]$. We derive a simple sufficient criterion for the unboundedness of the sum when $\urep{}$ is an arbitrary continuous unitary representation (cf. \thref{unbounded 1-cocycle}) and focus afterwards on the case of a finite tensor product representation (cf. \thref{case H^1 not zero}).\\ 

Let $G$ be a completely metrizable group and $(\pi,\SH)$ be a continuous unitary representation of ~$G$.
\lemma{\thlabel{unbounded 1-cocycle}
Let $(e_n)_{n\in\N}$ be an orthonormal sequence in $\SH$ such that, for all $n$, we have \[e_n \in \gen{\pi(g)v-v:\ \ g\in G,v\in\SH}{lin}.\] Further, let $(a_n)_{n\in\N}$ be a sequence in $\C$ for which the sum $\beta(g):= \sum_{n\in\N} a_n [\pi(g)e_n-e_n]$ converges in $\SH$ for every $g\in G$. Then $\beta:G\to\SH$ defines a (continuous) 1-cocycle which is a 1-coboundary \iof the sequence $(a_n)_{n\in\N}$ is square summable.
}

\pf{That $\beta$ is a continuous 1-cocycle follows from \thref{sum of 1-coboundaries}. If $\sum_{n\in\N} |a_n|^2 <\infty$, then we have \eqn{\label{sum is coboundary}\sum_{n\in\N} a_n [\pi(g)e_n-e_n] = \pi(g)v-v} for the vector $v:= \sum_{n\in\N}a_n e_n$. Conversely, assume that $\beta$ is a 1-coboundary. Then, we find some $v \in\SH$ which satisfies \eqref{sum is coboundary} for all $g\in G$. Let $(e_j)_{j\in J}$ be a complete orthonormal system in $\SH$ containing the $e_n$, i.e. we have $\N\subseteq J$. We expand $v$ w.r.t. the system $(e_j)$ and obtain the coefficients $b_j := \skp{v}{e_j}{}$. For $j\in J\backslash \N$, we put $a_j:=0$. Our assumption then leads to \[\sum_{j\in J}(a_j-b_j)[\pi(g)e_j-e_j] = 0 \quad\text{for all } g\in G.\] Therefore, we have $\sum_{j\in J}(a_j-b_j)\skp{e_j}{w}{} = 0$ for any \\$w\in D:=\gen{\pi(g)v-v:\ \ g\in G,v\in\SH}{lin}$. Since $e_n\in D$ by assumption, we obtain that \[a_n-b_n = \sum_{j\in J}(a_j-b_j)\skp{e_j}{e_n}{} = 0\] and we conclude that the sequence $(a_n)_{n\in\N}$ is square integrable. 
}

Now, we turn to  unbounded 1-cocycles of finite tensor products: Assume that we are given $m\geq 2$ continuous unitary representations $(\pi_i,\SH_i)_{i=1,\ldots,m}$ of $G$. We form the Hilbert tensor product $(\pi,\SH):= \widehat{\bigotimes}_{i=1,\ldots,m} (\pi_i,\SH_i)$.

\lemma{\thlabel{1-cocycle for tensor product}
Let $(e_n^{(i)})_{n\in\N}$ be an orthonormal sequence in $\SH_i$. For simplicity we will write $g.e_n^{(i)}$ instead of $\pi_i(g)e_n^{(i)}$. Assume that $(a_n)_{n\in\N}$ is a sequence in $\C$ such that, for all $i$ and all $g\in G$,
\[\sum_{n\in\N}|a_n|^2\mynorm{g.e_n^{(i)}-e_n^{(i)}}{\SH_i}^2 <\infty.\]
Then \[\beta(g):= \sum_{n\in\N}a_n [g.e_n^{(1)}\hotimes \ldots \hotimes g.e_n^{(m)}- e_n^{(1)}\hotimes \ldots \hotimes e_n^{(m)}]\]
defines a 1-cocycle w.r.t. the tensor product representation. Moreover, we have the estimate
\eqn{\label{estimate} \sum_{n\in\N}|a_n|^2\mynorm{g.e_n^{(1)}\hotimes \ldots \hotimes g.e_n^{(m)}- e_n^{(1)}\hotimes \ldots \hotimes e_n^{(m)}}{}^2 \leq m\cdot \sum_{i=1}^m \sum_{n\in\N}|a_n|^2\mynorm{g.e_n^{(i)}-e_n^{(i)}}{\SH_i}^2  .}
} 

\pf{
It is enough to verify that the sum converges for all $g\in G$ (\thref{sum of 1-coboundaries}). Choose some $g\in G$ and unitary operators $U_i \in U(\SH_i)$. For fixed $j\in \{1,\ldots,m\}$ and arbitrary $M>N \in \N$ we calculate
\eqnlines{rl}{ & \mynorm{\sum_{n=N}^M a_n [U_1e_n^{(1)}\hotimes\ldots\hotimes (g-\1).e_n^{(j)}\hotimes\ldots\hotimes U_me_n^{(m)}]}{}^2\\
= & \sum_{n,n' = N}^M a_n\oline{a_{n'}} \skp{U_1e_n^{(1)}}{U_1e_{n'}^{(1)}}{\SH_1}\ldots\skp{(g-\1).e_n^{(j)}}{(g-\1).e_{n'}^{(j)}}{\SH_j}\ldots\skp{U_me_n^{(m)}}{U_me_{n'}^{(m)}}{\SH_m} \\
= & \sum_{n=N}^M |a_n|^2\mynorm{(g-\1).e_n^{(j)}}{}^2 \underset{N,M\rightarrow\infty}{\longrightarrow} 0.}
We conclude that the sum $\sum_{n=1}^\infty a_n [U_1e_n^{(1)}\hotimes\ldots\hotimes (g-\1).e_n^{(j)}\hotimes\ldots\hotimes U_me_n^{(m)}]$ converges for every $j$ and $g$. This shows that 
\eqnlines{rl}{ & \sum_{n\in\N}a_n [g.e_n^{(1)}\hotimes \ldots \hotimes g.e_n^{(m)}- e_n^{(1)}\hotimes \ldots \hotimes e_n^{(m)}] \\ 
= & \sum_{j=1}^m \sum_{n\in\N} a_n [g.e_n^{(1)}\hotimes\ldots\hotimes g.e_n^{(j-1)}\hotimes (g-\1).e_n^{(j)}\hotimes e_n^{(j+1)}\hotimes\ldots\hotimes e_n^{(m)}]} converges for every $g\in G$.\\
The estimate \eqref{estimate} follows from
\eqnlines{rl}{ & \mynorm{g.e_n^{(1)}\hotimes \ldots \hotimes g.e_n^{(m)}- e_n^{(1)}\hotimes \ldots \hotimes e_n^{(m)}}{}^2 \\
= & \mynorm{\sum_{j=1}^m g.e_n^{(1)}\hotimes\ldots\hotimes g.e_n^{(j-1)}\hotimes (g-\1).e_n^{(j)}\hotimes e_n^{(j+1)}\hotimes\ldots\hotimes e_n^{(m)}}{}^2 \leq m \sum_{j=1}^m \mynorm{(g-\1).e_n^{(j)}}{\SH_j}^2.\qedhere }
}

For our purposes, tensor products of the form $(\pi_\mu^*,\SH_m^*)\hotimes(\pi_\lambda,\SH_\lambda)$ are of particular interest, where $(\pi_\lambda,\SH_\lambda)$ and $(\pi_\mu,\SH_\mu)$ are two irreducible unitary representations of $G$. The canonical isomorphism $\SH_\mu^*\hotimes\SH_\lambda \cong \SB_2(\SH_\mu,\SH_\lambda)$ induces a unitary representation on $\SB_2(\SH_\mu,\SH_\lambda)$ which is equivalent to the tensor product $\pi_\mu^*\hotimes \pi_\lambda$. It is given by the conjugation action \[(g,A)\mapsto \pi_\lambda(g)A\pi_\mu(g^{-1}),\text{ for all }g\in G, A\in \SB_2(\SH_\mu,\SH_\lambda). \] In particular, if $\urep{\lambda}$ and $\urep{\mu}$ are isomorphic unitary representations, then $\pi_\mu^*\hotimes \pi_\lambda$ is equivalent to the conjugation representation on $\SB_2(\SH_\lambda)$.

\lemma{\thlabel{unbounded 1-cocycles for conjugation}
Let $A\in \SB(\SH_\mu,\SH_\lambda)$ be a bounded linear operator such that \[\beta(g):= \pi_\lambda(g)A\pi_\mu(g^{-1}) - A \in \SB_2(\SH_\mu,\SH_\lambda)\] for every $g\in G$.
\enumeration{[i)]
\item Assume that $\urep{\lambda}$ and $\urep{\mu}$ are non-isomorphic representations. Then, $\beta$ is a 1-coboundary \iof $A\in \SB_2(\SH_\mu,\SH_\lambda)$.
\item Assume that $\urep{\lambda}=\urep{\mu}$. Then, $\beta$ is a 1-coboundary \iof $A\in \SB_2(\SH_\lambda) + \C\1$.
}}

\pf{In both cases, $\beta$ is a 1-coboundary \iof there exists $B\in \SB_2(\SH_\mu,\SH_\lambda)$ such that $\pi_\lambda(g)(A-B)\pi_\mu(g^{-1}) = A-B$ for all $g\in G$. We denote by $\SB^G(\SH_\mu,\SH_\lambda)$ the space of linear bounded operators intertwining the representations $\pi_\mu$ and $\pi_\lambda$. If the representation are nonisomorphic, then Schur's Lemma implies that $\SB^G(\SH_\mu,\SH_\lambda) = \{0\}$. This shows the first assertion. In the case $\pi_\mu = \pi_\lambda$, Schur's Lemma states that $\SB^G(\SH_\lambda) = \C\cdot\1$ which proves the second assertion.
}

\proposition{\thlabel{case H^1 not zero}
Let $\urep{\mu}$ and $\urep{\lambda}$ be two infinite dimensional, irreducible representations of the completely metrizable group $G$. Assume that there exists a bounded but not square summable sequence $(a_n)_{n\in\N}\in \ell^\infty(\N,\C)\backslash\ell^2(\N,\C)$ such that, for some orthonormal sequences $(e_n)_{n\in\N}$ and $(f_n)_{n\in\N}$ in $\SH_\lambda$ resp. $\SH_\mu$, 
\eqn{\label{trace class condition}\sum_{n\in\N}|a_n|^2\mynorm{(\pi_\lambda(g)-\1)e_n}{\SH_\lambda}^2 <\infty \quad\text{ and }\quad \sum_{n\in\N}|a_n|^2\mynorm{(\pi_\mu(g)-\1)f_n}{\SH_\mu}^2<\infty}
holds for all $g\in G$. Then $H^1(G,\pi_\mu^*\hotimes\pi_\lambda,\SH_\mu^*\hotimes\SH_\lambda) \neq \{0\}$.
}

\pf{
First, we note that $\mynorm{(\pi_\mu(g)-\1)f_n}{\SH_\mu} = \mynorm{(\pi_\mu^*(g)-\1)f^*_n}{\SH_\mu^*}$. By \thref{1-cocycle for tensor product}, we obtain a 1-cocycle $\beta(g):= \sum_{n\in\N}a_n [\pi_\mu^*(g)f_n^*\hotimes \pi_\lambda(g)e_n - f_n^*\hotimes e_n]$ for the tensor product. With respect to the conjugation action on $\SB_2(\SH_\mu,\SH_\lambda)$ this 1-cocycle has the form $\beta(g) = \pi_\lambda(g)A\pi_\mu(g^{-1}) - A$, where $A = \sum_{n\in\N}a_n e_n\otimes f_n^* = \sum_{n\in\N}a_n\skp{\cdot}{f_n}{\SH_\mu}e_n$. Since the sequence $(a_n)$ is bounded but not square summable, the operator $A$ is bounded but not Hilbert--Schmidt. If $\pi_\mu$ and $\pi_\lambda$ are not isomorphic, then \thref{unbounded 1-cocycles for conjugation} implies that $\beta$ is unbounded. If $\pi_\mu \cong \pi_\lambda$ are isomorphic, we may assume w.l.o.g. that $\pi_\mu=\pi_\lambda$. In this case, \thref{unbounded 1-cocycles for conjugation} implies that $\beta$ is unbounded unless $A \in \SB_2(\SH_\lambda)+\C\1$. If this is the case, then $(a_n)$ is the sum of a constant and a square summable sequence. It is clear that condition \eqref{trace class condition} then holds for all bounded sequences $(a_n)$ because the constant part is nonzero. In particular, we may choose the sequence $a_n:= (-1)^n$. The corresponding diagonal operator $A\in \SB(\SH_\lambda)$ cannot be written as a linear combination of the identity $\1$ and a Hilbert--Schmidt operator. In particular, the 1-cocycle $\beta$ is unbounded.
}

\section{Unitary highest weight representations of $U_p(\SH)$} 
\label{sec: U_p}

For finitely supported weights $\lambda \in\Z^\N$, the corresponding unitary highest weight representation $\urep{\lambda}$ of $U(\infty)$ extends to a norm-continuous unitary representation of $U_p(\SH)$ (for $p\in[1,\infty)$ and $\SH=\ell^2(\N,\C)$). We realize these representations as subrepresentations in finite tensor products of the identical representation and its dual representation of $U_p(\SH)$ on $\SH$ resp. on $\SH^*$ (cf. \thref{tensor product module}). Our results of the preceding section allow us to construct unbounded 1-cocycles for almost all finitely supported $\lambda$. There are only three exceptional cases where the first cohomology spaces are trivial (cf. ~\thref{H^1 for U_p}).\\ 

We denote by $\setlambda$ the set of all decreasingly ordered non-negative integer valued tuples with a finite number of positive entries.
To any $\lambda \in \setlambda$, we associate a \textit{Young diagram} $D_\lambda$ (also called \textit{Ferrers diagram}) which consists of $\ell_\lambda := \max\{j| \ \ \lambda_j >0\}$ rows and the $j$-th row has $\lambda_j$ boxes, so that the whole diagram consists of $|\lambda|:= \sum_{i=1}^{\ell_\lambda} \lambda_i$ boxes. Conversely, any Young diagram $D$ (with row (and column) length weakly decreasing) defines a unique $\lambda\in \setlambda$ by counting the row boxes. Given a Young diagram $D$, we obtain the \textit{conjugate} (or \textit{transposed}) Young diagram $D'$ by switching rows and columns. Thus, for any $\lambda\in \setlambda$, we define the \textit{conjugate} tuple $\lambda'\in\setlambda$ via the relation $D_{\lambda'} := D'_\lambda$. Note that $\ell_{\lambda'} = \lambda_1$ and the entry $j$ occurs $\lambda_j - \lambda_{j+1}$ times in $\lambda'$.

\example{For $\lambda := (3,2,2,1)$ we have 
\[D_\lambda = \ydiagram{3,2,2,1} \quad\text{ and }\quad D'_\lambda = \ydiagram{4,3,1}\]
and the conjugate tuple is given by $\lambda' =(4,3,1)$.
}

The corresponding \textit{Young tableau} $T_\lambda$ is obtained by filling in the boxes of the young diagram $D_\lambda$ with the numbers $1,2,\ldots,|\lambda|$ in the following manner: The number $1$ is placed in the top box of the first column. The number $k+1$ is placed in the box directly below $k$, if it exists and otherwise in the top box of the next column. 
\example{
For $\lambda = (3,2,2,1)$ this yields
\[
D_\lambda = \ydiagram{3,2,2,1} \quad\longrightarrow \quad
T_\lambda= \begin{ytableau}
1 & 5 & 8 \\
2 & 6 \\
3 & 7 \\
4 
\end{ytableau}
\]
}

Let $S_\lambda$ be the permutation group of the set $\{1,2,\ldots,|\lambda|\}$. Denote by $R_\lambda$ the subgroup of permutations leaving all subsets defined by the rows of $T_\lambda$ invariant and accordingly by $C_\lambda$ the subgroup leaving all subsets defined by the columns invariant. \\

Put $\SH:= \ell^2(\N,\C)$ and consider the $|\lambda|$-fold tensor product $\SH^{\hotimes |\lambda|}$. The permutation group $S_\lambda$ acts unitarily on the tensor product via \[(\sigma, \hotimes_{j=1}^{|\lambda|} v_j) \mapsto \hotimes_{j=1}^{|\lambda|} v_{\sigma^{-1}(j)}\] for $\sigma \in S_\lambda$ and $v_j \in \SH$. We denote this representation by $(\rho,\SH^{\hotimes |\lambda|})$. The group $U_p(\SH)$ (for $p\in [1,\infty]$) also acts on the tensor product space via  the tensor product representation $(\id,\SH)^{\hotimes |\lambda|}$. Both representations commute, i.e. for $\sigma \in S_\lambda$ and $g\in U_p(\infty)$, we have $\rho(\sigma)\circ \id^{\hotimes |\lambda|}(g) = \id^{\hotimes |\lambda|}(g) \circ \rho(\sigma)$. Therefore, the linear operator
\[P_\lambda := \sum_{r\in R_\lambda ,c\in C_\lambda} \sgn(c)\rho(cr),\] where $\sgn(c)$ denotes the signum of the permutation $c$,
commutes with $\id^{\hotimes |\lambda|}$ and its image
\[\SH_\lambda := P_\lambda (\SH^{\hotimes |\lambda|})\]
is a $U_p(\SH)$-invariant subspace. We write $\pi_\lambda$ for the corresponding continuous unitary representation on $\SH_\lambda$. Note that $P_\lambda$ is not necessarily an orthogonal projection, but the operators 
\[P_{C_\lambda} := \frac{1}{f_{C_\lambda}}\sum_{c\in C_\lambda} \sgn(c)\rho(c) \quad\text{ and }\quad P_{R_\lambda} := \frac{1}{f_{R_\lambda}}\sum_{r\in R_\lambda} \rho(r)\] with constants $f_{C_\lambda} := \prod_{j=1}^{\ell_{\lambda'}} (\lambda'_j!)$ and $f_{R_\lambda}:=\prod_{j=1}^{\ell_\lambda}(\lambda_j!) $, are orthogonal projections and we have $P_\lambda = \left(f_{C_\lambda}f_{R_\lambda}\right)\cdot P_{C_\lambda}P_{R_\lambda}$. Let $(e_n)_{n\in\N}$ denote an ONB for $\SH$. For $n\in\N$, we define the vectors 
\eqn{\label{e_k^lambda} e_n^{(\lambda)} := \sqrt{f_{C_\lambda}}\cdot\hotimes_{m=1}^{\ell_{\lambda'}} \wedge_{i=1}^{\lambda'_m} e_{n+i} = \sqrt{f_{C_\lambda}}\cdot P_{C_\lambda}(\hotimes_{m=1}^{\ell_{\lambda'}}\hotimes_{i=1}^{\lambda'_m} e_{n+i}).}
By construction, the vectors $e_n^{(\lambda)}$ are mutually orthgonal. Since $ P_{R_\lambda}(\hotimes_{m=1}^{\ell_{\lambda'}}\hotimes_{i=1}^{\lambda'_m} e_{n+i}) =\hotimes_{m=1}^{\ell_{\lambda'}}\hotimes_{i=1}^{\lambda'_m} e_{n+i} $, we have constructed an orthonormal sequence in $\SH_\lambda$. In particular, if $|\lambda|\geq 1$, the space $\SH_\lambda$ is infinite dimensional.

\example{
For $\lambda = (3,2,2,1)$ the orthonormal sequence $(e_n^{(\lambda)})_{n\in\N}$ in $\SH_\lambda$ is given by 
\[e_n^{(\lambda)} = \underbrace{\sqrt{4!3!1!}}_{=12}\cdot(e_{n+1}\wedge e_{n+2}\wedge e_{n+3}\wedge e_{n+4})\hotimes(e_{n+1}\wedge e_{n+2}\wedge e_{n+3})\hotimes e_{n+1}.\] 
}

\lemma{\thlabel{tableau-representation}
For any $\lambda\in\setlambda$, the continuous unitary representation $\urep{\lambda}$ of $U_p(\SH)$ extends the highest weight representation of $U(\infty)$ with highest weight $\lambda$ from \thref{HWR for U(infty)}. In particular, it is irreducible.
}

\pf{
For any $n\in\N$ with $n\geq \ell_\lambda$, consider the canonical embedding $\left(\C^n\right)^{\hotimes |\lambda|} \hookrightarrow \SH^{\hotimes |\lambda|}$. The operator $P_\lambda$ leaves the subspace $\left(\C^n\right)^{\hotimes |\lambda|}$ invariant as well as the restriction of $\id^{\hotimes |\lambda|}$ to the subgroup $U(n)\subset U_p(\SH)$. Therefore, the subspace $P_\lambda\left(\left(\C^n\right)^{\hotimes |\lambda|}\right)$ is $U(n)$-invariant and, according to Theorem ~9.3.9 in \cite{GW98}, this is a unitary highest weight module of $U(n)$ with highest weight $\lambda^{(n)}$. This shows that $\pi_\lambda|_{U(\infty)} = \indlim \pi_{\lambda^{(n)}}$ is the direct limit representation on $P_\lambda\Big(\left(\C^{(\N)}\right)^{\hotimes |\lambda|}\Big)$ which is dense in $\SH_\lambda$. }

Now, let $\lambda\in\Z^{(\N)}$ be an integer valued tuple such that all but finitely many entries are zero. For any such $\lambda\in\Z^{(\N)}$, we define $\lambda^\pm \in \setlambda$ to be those tuples for which the entries of the finitely supported tuples $\max(\pm \lambda,0)$ are decreasingly ordered. We put $|\lambda|:= |\lambda^+|+|\lambda^-| = \sum_{i=1}^\infty |\lambda_i|$. 

\definition{\thlabel{tensor product module} For any $\lambda\in\Z^{(\N)}$, we call the continuous unitary representation $\urep{\lambda} := \starurep{\lambda^-}\hotimes \urep{\lambda^+}$ \textit{unitary highest weight representation} of the group $U_p(\SH)$ \textit{with heighest weight $\lambda$}.}

We briefly comment on this definition.
\remark{That the representation $\urep{\lambda}$ is irreducible follows from Section 2.17 of \cite{Ol90}. Using Corollary I.14 of \cite{Ne98}, one finds that $\urep{\lambda}$ is a unitary highest weight representation of $U_p(\SH)$ in the sense of Definition III.6 in \cite{Ne98}. For the case $p=\infty$, we refer to \cite[Thm. 2.2]{BN12}.
}

In the remainder of this section we are going to show that, for $|\lambda|\geq 2$, the first cohomology space $H^1(U_p(\SH),\pi_\lambda,\SH_\lambda)$ never vanishes.

\lemma{\thlabel{U_p property}
Put $q:= \frac{2p}{p-2}$ if $2<p<\infty$ and $q:=\infty$ if $1\leq p\leq 2 $. For every sequence $a= (a_n)_{n\in\N}\in \ell^q(\N,\C)$, we have
\[\sum_{n\in\N} |a_n|^2\mynorm{(g-\1)e_n}{\SH}^2 <\infty\]
for every $g\in U_p(\SH)$. Moreover, this expression depends continuously on $g$ in the $U_p(\SH)$-topology.
}

\pf{
For $A\in\SB_r(\SH)$ and $r\in [1,\infty]$, we denote the $r$-th Schatten norm by \[\mynorm{A}{r}:= \case{\left(\Tr(|A|^r)\right)^{\frac{1}{r}} & \text{ if } r<\infty \\ \opnorm{A} &\text{ if } r=\infty}.\] We remind the reader of the generalized H\"older inequality (cf. \cite[Th. IV.11.2]{GGK00})
\eqn{\label{Hoelder} \mynorm{BC}{r} \leq \mynorm{B}{s}\cdot\mynorm{C}{t} \quad\text{for }B\in \SB_s(\SH), C\in \SB_t(\SH) \text{ and }\frac{1}{r}\leq \frac{1}{s} + \frac{1}{t}, \ \ s,t\in [1,\infty]. }

For $1\leq p\leq 2$, the statement is clear, since $g-\1\in\SB_2(\SH)$ and $\sum_{n\in\N} \mynorm{(g-\1)e_n}{}^2 = \mynorm{g-\1}{2}^2$. For $p>2$, we consider the diagonal operator $Ae_n = a_ne_n$ and note that $A\in \SB_q(\SH)$ and $\sum_{n\in\N} |a_n|^2\mynorm{(g-\1)e_n}{}^2 = \mynorm{(g-\1)A}{2}^2$. Applying \eqref{Hoelder} with $r=2$, $s=p$ and $t=q$ yields the assertion.
}

\proposition{\thlabel{H^1 not 0 for U_p}
Let $\lambda \in \Z^{(\N)}$ such that $|\lambda| \geq 2$. Then , for $p\in [1,\infty)$, we have \[H^1(U_p(\SH),\pi_\lambda,\SH_\lambda) \neq \{0\}.\]
}

\pf{
Put $a_n:= \frac{1}{\sqrt{n}}$. Then, the sequence $a=(a_n)_{n\in\N}$ is in $\ell^q(\N,\C)$ for all $q>2$ but not in $\ell^2(\N,\C)$ and we know from the previous lemma that the sum $\sum_{n\in\N}\frac{1}{n}\mynorm{(g-\1)e_n}{}^2$ is finite for each $g\in G$.
\aufzaehlung{
\item[$\text{The case } |\lambda^-|=0$:]
We first treat the case $\lambda \in\setlambda $. Applying \thref{1-cocycle for tensor product}, we obtain a 1-cocycle 
\[\widetilde{\beta}(g):= \sum_{n\in\N}\frac{1}{\sqrt{n}}[\hotimes_{m=1}^{\ell_{\lambda'}} \hotimes_{i=1}^{\lambda'_m} ge_{k+i} - \hotimes_{m=1}^{\ell_{\lambda'}} \hotimes_{i=1}^{\lambda'_m} e_{k+i}]\] for the tensor product representation $(\id,\SH)^{\hotimes |\lambda|}$. Projecting onto $\SH_\lambda$ and using \eqref{e_k^lambda} gives
\[\beta(g):= \sqrt{f_{C_\lambda}}\cdot P_{C_\lambda}(\widetilde{\beta}(g)) = \sum_{n\in\N} \frac{1}{\sqrt{n}}[g.e_n^{(\lambda)} - e_n^{(\lambda)}]. \] 
Here, we have used \eqref{e_k^lambda} and the fact that $P_{C_\lambda}$ commutes with the representation $\id^{\hotimes |\lambda|}$. Hence, we obtain a 1-cocycle for the highest weight representation $\urep{\lambda}$ and \thref{unbounded 1-cocycle} (applied to the sequence $(e_n^{(\lambda)})_{n\in\N}$) tells us that $\beta$ is in fact unbounded. The condition on the $e_n^{(\lambda)}$ is satisfied, since each $e_n^{(\lambda)}$ is a nontrivial eigenvector of all diagonal operators $t-\1$ for $t\in T_p(\SH)$, the maximal torus in $U_p(\SH)$.

\item[$\text{The case }|\lambda^+|=0$:] In this case $\urep{\lambda} = \starurep{\lambda^-}$ and the assertion follows from the fact that the 1-cohomology of the dual representation vanishes \iof it vanishes for the original representation.

\item[$\text{The remaining case } |\lambda^+|, |\lambda^-| \geq 1$:]. We want to apply \thref{case H^1 not zero}. It merely remains to verify that 
\[\sum_{n\in\N}\frac{1}{n}\mynorm{(g-\1).e_n^{(\lambda^\pm)}}{\SH_{\lambda^\pm}}^2 <\infty\]
for each $g\in U_p(\SH)$. Indeed, we can estimate for $\lambda\in\{\lambda^\pm\}$ as follows: 
\[\sum_{n\in\N}\frac{1}{n}\mynorm{(g-\1).e_n^{(\lambda)}}{\SH_{\lambda}}^2 = f_{C_\lambda}\sum_{n\in\N}\frac{1}{n}\mynorm{P_{C_\lambda}[\hotimes_{m=1}^{\ell_{\lambda'}} \hotimes_{i=1}^{\lambda'_m} ge_{k+i}-\hotimes_{m=1}^{\ell_{\lambda'}} \hotimes_{i=1}^{\lambda'_m} e_{k+i}]}{\SH_{\lambda}}^2\]
\[ \leq f_{C_\lambda}\sum_{n\in\N}\frac{1}{n}\mynorm{\hotimes_{m=1}^{\ell_{\lambda'}} \hotimes_{i=1}^{\lambda'_m} ge_{k+i}-\hotimes_{m=1}^{\ell_{\lambda'}} \hotimes_{i=1}^{\lambda'_m} e_{k+i}}{\SH_{\lambda}}^2 \overset{\eqref{estimate}}{\leq} f_{C_\lambda}|\lambda|^2 \sum_{n\in\N}\frac{1}{n}\mynorm{(g-\1)e_n}{\SH}^2 <\infty.\qedhere\]
}}

\thref{H^1 not 0 for U_p} together with \thref{sec: direct limit of compact groups}, \thref{identical representation} and \thref{dual representation} yields the following theorem.

\theorem{\thlabel{H^1 for U_p}
Put $\SH:= \ell^2(\N,\C)$. Let $\lambda \in \Z^{(\N)}$ be a finitely supported integer valued sequence and $\urep{\lambda}$ the corresponding unitary highest weight representation of $U_p(\SH)$. Then,
\[H^1(U_p(\SH),\pi_\lambda,\SH_\lambda) = \{0\} \Longleftrightarrow |\lambda| \leq 1. \]
}

One immediate consequence of \thref{H^1 for U_p} is that the group $U_p(\SH)$ does not have property (FH). Since property (T) implies property (FH) for arbitrary topological groups, this entails that $U_p(\SH)$ does not have property (T).   

\corollary{\thlabel{no (T) for U_p}
For $1\leq p < \infty$, the group $U_p(\SH)$ neither has property (FH) nor has property (T).
}

\appendix

\section{Infinite sums of 1-coboundaries}
\label{sec: infinite sums}
Let $(\pi,\SH)$ be a continuous unitary representation of the topological group $G$. In general, a 1-cocycle $\beta:G\to\SH$ need not be continuous. Here, we show that, if $\beta$ is the pointwise limit of a sequence of continuous 1-cocycles and $G$ is completely metrizable, then a Baire category argument shows that $\beta$ is continuous.

\definition{
Let $X$ be a topological Hausdorff space. A subset $A\subseteq X$ is called \textit{of first (Baire) category} if it can be written as a countable union of nowhere dense subsets. Otherwise, $A$ is called a subset \textit{of second (Baire) category}. The space $X$ is called \textit{Baire space}, if every nonempty open subset is of second category. A function $f:X\to X'$ with values in topological Hausdorff space $X'$ is called a \textit{function of first class} or a \textit{Baire one function} if it is the pointwise limit of a sequence $(f_n)_{n\in\N}$ of continuous functions $f_n:X\to X'$.  
}

\example{According to the Baire Category Theorem , every complete metric space $(X,d)$ is a Baire space (cf. \cite[Th\'eor\`eme 1 de $§5$]{Bour58}). In particular, the subset $A=X$ is of second category.}

\theorem{\thlabel{first class functions}
Let $(X,d)$, $(X',d')$ be metric spaces and $f:X\to X'$ be a first class function. Then the points of discontinuity of $f$ form a subset of first category.
}

\pf{see e.g. \cite[$§45.3$]{Hau62}.}

\proposition{\thlabel{pointwise limit}
Let $(\pi,\SH)$ be a continuous unitary representation of the completely metrizable group $G$ and let $\beta:G\to\SH$ be the pointwise limit of a sequence of continuous 1-cocycles $(\beta_n)_{n\in\N}$. Then $\beta$ is a continuous 1-cocycles.
}

\pf{
The pointwise convergence implies that the map $\beta$ is a first class function. It is straightforward to check that the 1-cocycle relation \eqref{1-cocycle} is satisfied, hence $\beta$ is a 1-cocycle. Moreover, the 1-cocycle relation shows that the points of discontinuity of $\beta$ are either the empty set or coincide with the whole group $G$. In the latter case, the group $G$ would have to be of first category since $\beta$ is a first class function (\thref{first class functions}). This contradicts the fact that $G$ is of second category (Baire Category Theorem).
}

\corollary{\thlabel{sum of 1-coboundaries}
Let $(\pi,\SH)$ be a continuous unitary representation of the completely metrizable group $G$ and let $(v_n)_{n\in\N}$ be a sequence in $\SH$ for which the sum $\sum_{n\in\N}\pi(g)v_n-v_n$ converges for every $g\in G$. Then $\beta(g):= \sum_{n\in\N}\pi(g)v_n-v_n$ defines a continuous 1-cocycle.
}


\end{document}